\pdfoutput=1
\documentclass[11pt, reqno]{amsart}
\usepackage[dvips]{graphicx}
\usepackage{bmpsize} 
\usepackage[margin=1.12in]{geometry}

\usepackage{mathrsfs}
\usepackage{amsmath,amscd}
\usepackage{amssymb}
\usepackage{amscd}
\graphicspath{/home/jhois/Documents/Images/}  

\newtheorem{theorem}{Theorem}

\newtheorem{proposition}[theorem]{Proposition}
\newtheorem{lemma}[theorem]{Lemma}

\theoremstyle{definition} 

\newcommand{\R}{\mathbb R}

\newcommand{\Gcal}{\mathcal{G}}
\newcommand{\Kcal}{\mathcal{K}}

\begin{document}

\author[J.A.~Hoisington]{Joseph Ansel Hoisington} \address{Rose-Hulman Institute of Technology}\email{hoisingt@rose-hulman.edu}

\title[Spherical Images]{A Note on Spherical Images of Convex Hypersurfaces}


\date{}

\begin{abstract}
A theorem of Wu states that the closure of the image of the Gauss map of an open convex hypersurface is geodesically convex, although the image itself may not be.  We show that the general version of this theorem follows from a short proof for smooth hypersurfaces due to do Carmo and Lawson.  We also show that all closed, geodesically convex subsets of spheres are attained as spherical images of convex hypersurfaces and give examples of convex surfaces whose spherical images cannot be attained as spherical images of smooth surfaces. 	
\end{abstract}

\maketitle




Let $M^n$ be a connected, noncompact, convex hypersurface in $\R^{n+1}$.  In \cite{Wu}, Wu disproved the assumption, going back at least to work of Busemann \cite[p.25]{Bu1}, that the image of the Gauss map of $M$ must be geodesically convex in the unit sphere $S^n$.  He then proved: 
\begin{theorem}[Wu, \cite{Wu}]
\label{Wu_theorem}
Let $M$ be a connected, noncompact, convex hypersurface in $\R^{n+1}$ with Gauss map $\varphi:M \to S^n$ and $\Gcal$ the image of $\varphi$.  Then the closure $\overline{\Gcal}$ of $\Gcal$ is geodesically convex in $S^n$, in that it contains a minimizing geodesic segment between any two of its points. 
\end{theorem}
\noindent Wu established Theorem \ref{Wu_theorem} for all convex hypersurfaces, without assuming any regularity beyond that which follows from convexity, for a generally multivalued version of the Gauss map that takes $x \in M$ to the set of all outer unit normals to supporting hyperplanes at $x$.  In \cite{dCL}, do Carmo and Lawson gave a short proof of Theorem \ref{Wu_theorem} for smooth convex hypersurfaces.  Note that Wu gave examples of smooth convex surfaces in $\R^3$ whose Gauss maps have nonconvex image, so Theorem \ref{Wu_theorem} is in a sense the strongest result of its type among smooth hypersurfaces and in all dimensions $(n+1) \geq 3$.  In this note, we derive Theorem \ref{Wu_theorem} in full generality from do Carmo and Lawson's proof in the smooth case.  The basic idea of the derivation is that do Carmo and Lawson's argument applies to $C^{1,1}$ hypersurfaces, and therefore to the outer parallel $M_r$ to any convex hypersurface $M$ at any positive distance $r$, and that the Gauss map of $M_r$ has the same image as that of $M$.  For completeness, we present the proof in full below.  This proof makes clear that a version of Theorem \ref{Wu_theorem} also applies to unbounded convex sets with empty interior, for which the domain of the Gauss map is the set's normal cycle as defined in Z\"ahle's work \cite{Za1}. 
\smallskip 

Theorem \ref{Wu_theorem} would follow immediately from do Carmo and Lawson's work if the image of the Gauss map of every noncompact convex hypersurface could be attained as the image of the Gauss map of a smooth convex hypersurface, but this is not the case.  In Proposition \ref{non_smooth_prop}, we describe a family of convex surfaces in $\R^3$ whose Gauss map images are not the image of the Gauss map of any smooth surface.  The proof of Proposition \ref{non_smooth_prop} is based on the existence of convex plane domains $\Kcal$ whose boundary $\partial \Kcal$ is nondifferentiable at a dense set of points.  We refer to \cite[Example 2.8]{AB1} and \cite[Ch.2, Exercise 16]{HW1} for examples of such domains.  Note that convex curves are rectifiable, and thus differentiable on a set of full $1$-dimensional Hausdorff measure, and that a convex plane domain $\Kcal$ has a unique support line through any point $x$ at which its boundary $\partial \Kcal$ is differentiable.  Conversely, an elementary argument shows that the uniqueness of the support line to $\Kcal$ at $x$ implies $\partial \Kcal$ is differentiable there.    

\begin{lemma}
\label{non_smooth_lemma}
A compact, convex set $\Kcal$ in $\R^2$ whose boundary $\partial \Kcal$ is nondifferentiable on a dense set is not the image of a smooth map of a second-countable differentiable manifold. 
\end{lemma}

\begin{proof} Let $\Kcal$ be such a convex set, and suppose there was a smooth map $F$ of a second-countable differentiable manifold $Y$ with image $\Kcal$.  For $x \in \partial \Kcal$ at which $\partial \Kcal$ is not differentiable, as discussed above, there are two distinct support lines to $\Kcal$ through $x$, corresponding to linearly independent functionals $h_{\vec{u}_1} = \langle \vec{u}_1, \ast \rangle$ and $h_{\vec{u}_2} = \langle \vec{u}_2, \ast \rangle$ on $\R^2$ whose maximum on $\Kcal$ is attained uniquely at $x$.  For $y \in F^{-1}(x)$, both functions $h_{\vec{u}_i} \circ F$ would then attain their maximum on $Y$ at $y$, which would imply that $\langle dF_y (\vec{w}), \vec{u}_i \rangle = 0$ for all tangent vectors $\vec{w}$ to $Y$ at $y$, and therefore, because $\vec{u}_1$, $\vec{u}_2$ are linearly independent, that $dF_y = 0$.  
\smallskip 

Given a countable compact exhaustion $\mathscr{C}_1 \subseteq \mathscr{C}_2 \subseteq \mathscr{C}_3 \subseteq \cdots$ of $Y$, the Baire category theorem would imply that $F(\mathscr{C}_{m}) \cap \partial \Kcal$ has nonempty interior in $\partial \Kcal$ for sufficiently large $m$. For a nonempty open interval $\mathcal{I}$ in such an $F(\mathscr{C}_{m}) \cap \partial \Kcal$, the density of nondifferentiable points of $\partial \Kcal$ in $\mathcal{I}$ would imply that for all $x \in \mathcal{I}$, there would be $y_0 \in F^{-1}(x) \cap \mathscr{C}_m$ that is the limit of a sequence $(y_j)_{j=1}^{\infty}$ in $\mathscr{C}_m$ such that  $\partial \Kcal$ is nondifferentiable at $F(y_j)$, and therefore, as shown above, $dF_{y_j} = 0$.  By continuity, we would then have $dF_{y_0} = 0$.  Let $l$ be a line in $\R^2$ such that the image of $\mathcal{I}$ via orthogonal projection $\pi_l$ onto $l$ has nonempty interior---for example, a support line to $\Kcal$ through a point in $\mathcal{I}$.  Because all $x \in \mathcal{I}$ would have $y_0 \in F^{-1}(x)$ with $dF_{y_0} = 0$, all $x \in \mathcal{I}$ would map via $\pi_l$ to critical values of $\pi_l \circ F$, contradicting Sard's theorem.  $\Kcal$ therefore cannot be the image of such a smooth map. \end{proof} 

\begin{lemma}
\label{Gauss_image_lemma}
A closed subset of the unit sphere $S^n$ in $\R^{n+1}$ is the image of the Gauss map of a connected, noncompact, convex hypersurface if and only if it is geodesically convex and contained in a closed hemisphere.  
\end{lemma}

\begin{proof} Theorem \ref{Wu_theorem} implies that a closed spherical image of a connected convex hypersurface is geodesically convex, and because the image of the Gauss map of a noncompact convex hypersurface is not the whole of $S^n$, it must be contained in a closed hemisphere, cf. \cite[Theorem 1]{Wu}.  To establish the converse, let $\widehat{\Kcal}$ be a closed, geodesically convex subset of a closed hemisphere.  For $\vec{u} \in S^n$, let $H_{\vec{u}} = \lbrace p \in \R^{n+1} \ : \ \langle p, \vec{u} \rangle \leq 1 \rbrace$, and define $\Omega$ to be: 
\begin{align*}
\displaystyle \bigcap\limits_{\vec{u} \in \widehat{\Kcal}} H_{\vec{u}}. 
\end{align*}
\noindent Let $M = \partial \Omega$, and let $\varphi:M \to S^n$ be the Gauss map of $M$.  It is immediate from the construction of $\Omega$ that for all $\vec{u} \in \widehat{\Kcal}$, we have $\vec{u} \in M$ and $\varphi(\vec{u}) = \vec{u}$, so that $\widehat{\Kcal} \subseteq Im(\varphi)$.  For $\vec{u} \in S^n \setminus \widehat{\Kcal}$, on the other hand, $h_{\vec{u}} = \langle \vec{u}, \ast \rangle$ is unbounded on $\Omega$, which implies that $\vec{u} \notin Im(\varphi)$, and therefore that $\widehat{\Kcal} = Im(\varphi)$.  The fact that $h_{\vec{u}}$ is unbounded on $\Omega$ is immediate for those $\vec{u}$ whose distance from $\widehat{\Kcal}$ in $S^n$ is at least $\frac{\pi}{2}$.  For those $\vec{u} \in S^n \setminus \widehat{\Kcal}$ whose distance from $\widehat{\Kcal}$ is less than $\frac{\pi}{2}$, the compactness and geodesic convexity of $\widehat{\Kcal}$ imply by an elementary argument that there is a unique $\vec{v} \in \widehat{\Kcal}$ nearest to $\vec{u}$.  The $2$-plane spanned by $\vec{u}$ and $\vec{v}$ in $\R^{n+1}$ intersects $M$ in an unbounded convex curve that includes a ray based at $\vec{v}$ and orthogonal to it, and $h_{\vec{u}}$ is unbounded on this ray, and therefore on $\Omega$. \end{proof}

\begin{proposition}
\label{non_smooth_prop} 
There are convex surfaces in $\R^3$ whose Gauss map images are not the image of the Gauss map of any smooth convex surface. 
\end{proposition}

\begin{proof} Let $S^2$ be the unit sphere in $\R^3$, and define $f:\R^2 \to S^2$ as follows: for $p \in \R^2$, letting $(p,1)$ be the point in $\R^3$ with the same $x$- and $y$-coordinates as $p$ and $z$-coordinate $1$, 
\begin{align*}
\displaystyle f(p) = \frac{(p,1)}{\sqrt{1 + |p|^2}}. 
\end{align*}
\noindent Note that $f$ is a diffeomorphism from $\R^2$ to the open hemisphere consisting of those $\vec{u} \in S^2$ with positive $z$-coordinate.  In synthetic terms, $f$ is essentially the embedding of $\R^2$ into its projective completion.  In particular, $f$ maps convex subsets of $\R^2$ to geodesically convex subsets of $S^2$.  Let $\Kcal$ be a compact, convex subset of $\R^2$ whose boundary has a dense set of nondifferentiable points, as in Lemma \ref{non_smooth_lemma}, and let $\widehat{\Kcal} = f(\Kcal)$.  Lemma \ref{Gauss_image_lemma} implies there is a noncompact convex hypersurface $M$ in $\R^3$ whose Gauss map has image $\widehat{\Kcal}$.  If $\widetilde{M}$ were a smooth convex surface whose Gauss map $\widetilde{\varphi}$ had image $\widehat{\Kcal}$, then $f^{-1} \circ \widetilde{\varphi}$ would give a smooth map of $\widetilde{M}$ with image $\Kcal$, contradicting Lemma \ref{non_smooth_lemma}, so no such surface exists. \end{proof}

\noindent Note that the existence of $C^{1,1}$ convex surfaces whose Gauss map image cannot be attained by a smooth surface follows from the existence of any such convex surface by the outer parallel construction in the proof of Theorem \ref{Wu_theorem} below, but that the surfaces in the proof of Lemma \ref{Gauss_image_lemma} are $C^{1,1}$ by construction.  The proofs of Lemma \ref{Gauss_image_lemma} and Proposition \ref{non_smooth_prop} therefore give explicit examples of $C^{1,1}$ surfaces with this property. 

\begin{proof}[Proof of Theorem \ref{Wu_theorem}] We first prove the theorem for $C^{1,1}$ hypersurfaces, following do Carmo and Lawson's proof in the smooth case.  Let $\vec{u}_0,\vec{u}_1 \in \Gcal$, suppose $\vec{u}_0$ and $\vec{u}_1$ are not antipodal, and let $\vec{u}$ be a point on the minimizing great circle arc between them.  Let $h_{\vec{u}} = \langle \vec{u}, \ast \rangle$, and let $\nabla h_{\vec{u}}$ be the gradient of $h_{\vec{u}}|_{M}$; that is, at $x \in M$,  
\begin{align}
\label{pf_eqn_1}
\displaystyle \nabla h_{\vec{u}} = \vec{u} - \langle \varphi(x), \vec{u} \rangle \varphi(x)
\end{align}	
Because $M$ is $C^{1,1}$, its Gauss map is Lipschitz, and it follows from the identity in \eqref{pf_eqn_1} that $\nabla h_{\vec{u}}$ is a Lipschitz section of the tangent bundle of $M$.  Because $\nabla h_{\vec{u}}$ is Lipschitz, it has unique trajectories in a neighborhood of any $x \in M$, as discussed in \cite[Chapter 1]{CL1}.  Note that solutions of ordinary differential equations with coefficients that are continuous but not Lipschitz are generally not unique, cf. \cite[p.7-8]{CL1}, so the uniqueness of the trajectories of $\nabla h_{\vec{u}}$ would not follow in this manner for a hypersurface that was $C^{1}$ but not $C^{1,1}$.  With the local existence and uniqueness of trajectories of $\nabla h_{\vec{u}}$ established, the proof in the $C^{1,1}$ case proceeds as in do Carmo and Lawson's work: let $\psi_{\vec{u}}(x,t)$ be the trajectory of $\nabla h_{\vec{u}}$ with $\psi_{\vec{u}}(x,0) = x$.  Because $|\nabla h_{\vec{u}}| \leq 1$ and $M$ is complete, $\psi_{\vec{u}}(x,t)$ exists for all $t \in \R$, and
\begin{align}
\displaystyle h_{\vec{u}} \circ \psi_{\vec{u}}(x,t) - h_{\vec{u}} \circ \psi_{\vec{u}}(x,0) = \int_{0}^{t} |\nabla h_{\vec{u}} \circ \psi_{\vec{u}}(x,\tau)|^{2} d\tau. 
\end{align}
If $|\nabla h_{\vec{u}} \circ \psi_{\vec{u}}(x,t)|$ were bounded away from $0$ for $t > 0$, we would therefore have: 
\begin{align}
\displaystyle \lim\limits_{t \to \infty} h_{\vec{u}} \circ \psi_{\vec{u}}(x,t) = \infty.  
\end{align}
The supporting hyperplanes to $M$ normal to $\vec{u}_0$ and $\vec{u}_1$ intersect, however, since $\vec{u}_0$ and $\vec{u}_1$ are not antipodal, which implies that $h_{\vec{u}}$ is bounded above on $M$.  There is therefore a sequence of points $(x_k)_{k = 1}^{\infty}$ along $\psi_{\vec{u}}(x,t)$ such that $|\nabla h_{\vec{u}}(x_k)|$ goes to $0$, and therefore $\varphi(x_k) \to \vec{u}$, which implies that $\vec{u}$, and therefore the minimizing great circle arc between $\vec{u}_0$ and $\vec{u}_1$, belongs to $\overline{\Gcal}$.  Having established that minimizing great circle arcs between nonantipodal points of $\Gcal$ belong to $\overline{\Gcal}$, the same conclusion follows for nonantipodal points of $\overline{\Gcal}$ by sequential approximation and for antipodal points $\pm \vec{u}_0$ of $\overline{\Gcal}$ by choosing a point $\vec{u}_1 \in \overline{\Gcal}$ distinct from both, and therefore not antipodal to either, which exists because $M$ is connected.  This completes the proof in the $C^{1,1}$ case.  For a general convex hypersurface $M$, any $r > 0$, all $x \in M$, and all outer unit normals $\vec{u}$ to $M$ at $x$, $x + r\vec{u}$ belongs to the outer parallel hypersurface $M_r$ to $M$ at distance $r$, which is $C^{1,1}$, cf. \cite[Theorem 1.2]{GH1}.  Conversely, all points in $M_r$ are of this form for a unique such $x \in M$ and $\vec{u}$.  Letting $\varphi_r$ be the Gauss map of $M_r$, the identity $\varphi_r(x + r\vec{u}) = \vec{u}$ implies that $\varphi_r$ has the same image as $\varphi$ and completes the proof in general. \end{proof}

\noindent As noted above, $C^{1,1}$ regularity is perhaps the weakest natural condition that guarantees the local existence and uniqueness of trajectories of the vector field in \eqref{pf_eqn_1}.  In this sense, it is the weakest condition that is compatible with do Carmo and Lawson's argument, which is based on the flow along these trajectories, without a more substantial modification.  It is also the strongest regularity condition that holds for outer parallels to arbitrary convex hypersurfaces. 
\medskip 

\noindent {\bf Acknowledgments and AI Use:} I am happy to thank Matteo Raffaelli for his feedback about this work, and for describing examples that established the conclusion of Proposition \ref{non_smooth_prop}.  These examples were found using standard AI tools, and the proof of Lemma \ref{non_smooth_lemma} is derived from the construction of these examples.  No other use of AI was made in this work. 




\begin{thebibliography}{ABCDE}
	
	
	\bibitem[AB96]{AB1} Stephanie B. Alexander and Richard L. Bishop: {\em Comparison theorems for curves of bounded geodesic curvature in metric spaces of curvature bounded above}, Differential Geometry and its Applications 6.1 (1996): 67-86.
	
	\bibitem[Bu58]{Bu1} Herbert Busemann: {\em Convex Surfaces}, Interscience, New York, 1958
	
	\bibitem[CL55]{CL1} Earl A. Coddington and Norman Levinson: {\em Theory of Ordinary Differential Equations}, McGraw-Hill (1955).
	
	\bibitem[dCL72]{dCL} Manfredo do Carmo and H. Blaine Lawson: {\em Spherical Images of Convex Surfaces}, Proc. Amer. Math Soc. 31 (1972): 635-636. 

	
	\bibitem[GH14]{GH1} Mohammad Ghomi and Ralph Howard: {\em Tangent Cones and Regularity of Real Hypersurfaces}, J. Reine Angew. Math. 697 (2014), 221–247.	
	
	\bibitem[HW20]{HW1} Daniel Hug and Wolfgang Weil: {\em Lectures on Convex Geometry}, Vol. 286. Cham: Springer, 2020.
	
	
	
	\bibitem[Wu74]{Wu} H. Wu: {\em The Spherical Images of Convex Hypersurfaces}, Journal of Differential Geometry 9.2 (1974): 279-290.
	
	\bibitem[Z\"a87]{Za1} Martina Z\"ahle: {\em Curvatures and Currents for Unions of Sets with Positive Reach}, Geometriae Dedicata 23.2 (1987): 155-171.
\end{thebibliography}
\end{document}